\documentclass[a4paper,12pt, reqno]{amsart} 

\usepackage{latexsym}
\usepackage{amssymb,amsfonts,amsmath,mathrsfs}

\newcommand{\ep}{\varepsilon}
\newtheorem{theorem}{Theorem}
\newtheorem{corollary}{Corollary}
\newtheorem{lemma}{Lemma}

\begin{document}
\title{On simple zeros of the Riemann zeta-function}
\author{H. M. Bui \and D. R. Heath-Brown}
\address{Institut f\"ur Mathematik, Universit\"at Z\"urich, Z\"urich
  CH-8057, Switzerland} 
\email{hung.bui@math.uzh.ch}
\address{Mathematical Institute, University of Oxford, Oxford OX1 3LB,
  United Kingdom} 
\email{rhb@maths.ox.ac.uk}

\begin{abstract}
We show that  at least $19/27$ of the zeros of the
Riemann zeta-function are simple, assuming the Riemann Hypothesis
(RH). This was previously established by Conrey, Ghosh
and Gonek [Proc. London Math. Soc. \textbf{76} (1998), 497--522] under
the additional assumption of the Generalised
Lindel\"of Hypothesis (GLH).  We are able to remove this hypothesis by
careful use of the generalised Vaughan identity.
\end{abstract}
\maketitle

\section{Introduction}

An important question in number theory is to understand the
distribution of the zeros of the Riemann zeta-function. In this paper,
we study the simple zeros on the critical line.  

Let $N(T)$ denote the number of zeros $\rho=\beta+i\gamma$ with
$0<\gamma<T$, where each zero is counted with multiplicity, denoted by
$m(\rho)$. Let $N^{*}(T)$ denote the number of such zeros which are
simple ($m(\rho)=1$), and let $N_d(T)$ denote the number of such
distinct zeros (i.e. each zero is counted precisely once without
regard to its multiplicity). Define $\kappa^{*}$ and $\kappa_d$ by 
\begin{equation*}
\kappa^{*}:=\liminf_{T\rightarrow\infty}\frac{N^{*}(T)}{N(T)},\qquad
\kappa_d:=\liminf_{T\rightarrow\infty}\frac{N_d(T)}{N(T)}. 
\end{equation*}

Unconditionally, it is known that $\kappa^{*}\geq0.4058$ (see
[\textbf{\ref{L}}, \textbf{\ref{H-B2}}, \textbf{\ref{C2}},
\textbf{\ref{A}}, \textbf{\ref{C1}}, \textbf{\ref{BCY}}]   
for results in this direction). Conditionally, using the pair
correlation of the zeros of the Riemann zeta-function Montgomery
[\textbf{\ref{M}}] showed that $\kappa^{*}\geq2/3$ and
$\kappa_d\geq5/6$ on RH. This was later improved by Cheer and Goldston
[\textbf{\ref{CG}}] to $\kappa^{*}\geq0.6727$ under the same
condition. Assuming RH and GLH, Conrey, Ghosh and Gonek
[\textbf{\ref{CGG1}}] showed that $\kappa^{*}\geq19/27$ and
$\kappa_d\geq0.84568$. Their paper used the mollifier method
(described in the next section). We also note that Montgomery's pair
correlation conjecture implies that almost all the zeros are simple.

In this paper, we use Heath-Brown's generalisation of the Vaughan
identity [\textbf{\ref{H-B1}}]
(see Lemma 3 in Section 3 below) to remove the GLH assumption in
the paper of Conrey, Ghosh and Gonek. As a result, we obtain 

\begin{theorem}
Assuming RH we have
\begin{equation*}
\kappa^{*}\geq\frac{19}{27}.
\end{equation*}
\end{theorem}

\begin{corollary}
Assuming RH we have
\begin{equation*}
\kappa_d\geq0.84665.
\end{equation*}
\end{corollary}

The corollary is a consequence of our theorem following an observation of
Montgomery [\textbf{\ref{M}}] that 
\begin{equation*}
2N^{*}(T)\leq\sum_{0<\gamma\leq
  T}\frac{(m(\rho)-2)(m(\rho)-3)}{m(\rho)}\leq\sum_{0<\gamma\leq
  T}m(\rho)-5N(T)+6N_d(T). 
\end{equation*}
Cheer and Goldston [\textbf{\ref{CG}}] also showed that
\begin{equation*}
\sum_{0<\gamma\leq T}m(\rho)\leq\big(1.3275+o(1)\big)N(T).
\end{equation*}
Hence
\begin{equation*}
\kappa_d\geq\frac{5+2\kappa^{*}-1.3275}{6}\geq0.84665.
\end{equation*}

Before embarking on the proof we record one piece of notation that we
will use throughout the paper, namely that we will write $\mathscr{L}=\log T/2\pi$. 

\section{The setup}

To get a lower bound for $\kappa^{*}$, it suffices to consider the
first and second mollified moments of the derivative of the Riemann
zeta-function. This section is mostly a summary of
[\textbf{\ref{CGG1}}]. 

We first note that $\rho$ is a simple zero if and only if
$\zeta'(\rho)\ne0$. Hence it follows from Cauchy's inequality that 
\begin{equation}\label{1}
N^{*}(T)\geq\frac{\big|\sum_{0<\gamma\leq
    T}B\zeta'(\rho)\big|^2}{\sum_{0<\gamma\leq
    T}\big|B\zeta'(\rho)\big|^2}, 
\end{equation}
for any regular function $B(s)$. Here we shall take $B(s)$ to be a
mollifier of the form 
\begin{eqnarray*}
B(s)=\sum_{k\leq y}\frac{b(k)}{k^s},
\end{eqnarray*}
where
\begin{equation}\label{bdef}
b(k)=\mu(k)P\bigg(\frac{\log y/k}{\log y}\bigg),
\end{equation}
with $P(x)$ being a polynomial with real coefficients satisfying
$P(0)=0$, $P(1)=1$, and $y=T^\vartheta$, $0<\vartheta<1/2$. 

The following result is essentially in [\textbf{\ref{CGG1}}]
\big(see (3.13), (3.21), (3.26), (3.27), (5.1), (5.4), (5.5)\big). 

\begin{lemma}
For any fixed $\ep>0$ we have
\begin{eqnarray*}
S_1&:=&\sum_{0<\gamma\leq T}B\zeta'(\rho)\\
&=&\frac{T\mathscr{L}^2}{2\pi}-\overline{\mathcal{M}_1}+
O(T\mathscr{L})+O_\varepsilon(yT^{1/2+\varepsilon})  
\end{eqnarray*}
and
\begin{eqnarray*}
S_2&:=&\sum_{0<\gamma\leq T}B\zeta'(\rho)B\zeta'(1-\rho)\\
&=&\frac{T\mathscr{L}^3}{2\pi}\bigg(\tfrac{1}{2}+
3\vartheta\int_{0}^{1}P(u)^2du\bigg)-2\emph{Re}(\mathcal{M}_2)+
O_\varepsilon(T\mathscr{L}^{2+\varepsilon})+  
O_\varepsilon(yT^{1/2+\varepsilon}), 
\end{eqnarray*}
where
\begin{equation}\label{2}
\mathcal{M}_\nu=\sum_{k\leq y}\sum_{m\leq
  kT/2\pi}\frac{a_\nu(m)b(k)}{k}e\bigg(\!-\frac{m}{k}\bigg). 
\end{equation}
Here $\nu=1$ or $2$, and the coefficients $a_\nu(m)$ are defined by
\begin{eqnarray}\label{3}
\frac{\zeta'}{\zeta}(s)\zeta'(s)=\sum_{n=1}^{\infty}\frac{a_1(n)}{n^s},
\qquad\frac{\zeta'}{\zeta}(s)\zeta'(s)^2B(s)=\sum_{n=1}^{\infty}\frac{a_2(n)}{n^s}. 
\end{eqnarray}
\end{lemma}

The main difficulty in the paper of Conrey, Ghosh and Gonek is to extract
the main terms and estimate the error terms in $\mathcal{M}_\nu$. At this point,
if we assume the Generalised Riemann Hypothesis (GRH), we can apply
Perron's formula to the sum over $m$ in \eqref{2}, and then move the
line of integration to $\textrm{Re}(s)=1/2+\varepsilon$. The main
terms arise from the residues of the pole at $s=1$ and the error
terms in this case are easy to handle. To avoid assuming GRH, however,
we first need to express the additive character $e(-m/k)$ in \eqref{2}
in terms of multiplicative characters, and then write $\mathcal{M}_\nu$ in the
following form \big(see [\textbf{\ref{CGG1}}; (5.12) and (5.14)]\big) 
\begin{eqnarray}\label{9}
\mathcal{M}_\nu=\sum_{q\leq y}\sum_{\psi(\textrm{mod}\
  q)}{\!\!\!\!\!\!}^{\textstyle{*}}\
\tau(\overline{\psi})\sum_{k\leq
  y/q}\frac{b(kq)}{kq}\sum_{d|kq}\delta(q,kq,d,\psi)\sum_{m\leq
  kqT/2\pi d}a_\nu(md)\psi(m), 
\end{eqnarray}
where $\sum^{*}$ denotes summation over all primitive characters
$\psi$(mod $q$), $\tau(\psi)$ is the Gauss sum, and  
\begin{equation}\label{deldef}
\delta(q,kq,d,\psi)=\sum_{l|(d,k)}\frac{\mu(d/l)}{\varphi(kq/l)}
\overline{\psi}\bigg(\frac{-k}{l}\bigg)
\psi\bigg(\frac{d}{l}\bigg)\mu\bigg(\frac{k}{l}\bigg). 
\end{equation}

Following [\textbf{\ref{CGG1}}] we choose a large constant $A$ and set
$\eta=\mathscr{L}^A$. We then split the $q$-summation into three
cases: $q=1$, $1<q\le\eta$, and $\eta<q\le y$.  We write
$\mathcal{M}_\nu=\mathcal{M}_{\nu,1}+\mathcal{M}_{\nu,2}+\mathcal{M}_{\nu,3}$
accordingly. The case 
$q=1$ gives rise to the main terms 
\big(see [\textbf{\ref{CGG1}}; Section 8]\big) 
\begin{equation}\label{4}
\mathcal{M}_{1,1}=\frac{T\mathscr{L}^2}{2\pi}
\bigg(\tfrac{1}{2}-\vartheta\int_{0}^{1}P(u)du\bigg)+O(T\mathscr{L}) 
\end{equation}
and
\begin{eqnarray}
\mathcal{M}_{2,1}&=&\frac{T\mathscr{L}^3}{2\pi}\bigg(\tfrac{1}{12}-
\tfrac{\vartheta}{2}\int_{0}^{1}P(u)du+\tfrac{3\vartheta}{2}
\int_{0}^{1}P(u)^2du\nonumber\\ 
&&\qquad\ \qquad-\tfrac{\vartheta^2}{2}
\bigg(\int_{0}^{1}P(u)du\bigg)^2-\tfrac{1}{24\vartheta}
\int_{0}^{1}P'(u)^2du\bigg)+O(T\mathscr{L}^2).\label{5}
\end{eqnarray}
The terms with $1<q\leq\eta$ are handled 
using Siegel's theorem  on exceptional real zeros of $L$-functions,
\big(see 
[\textbf{\ref{CGG1}}; (5.15) and (6.14)]\big) to give 
\begin{eqnarray}\label{6}
\mathcal{M}_{\nu,2}\ll_A T\exp\big(-c(A)\sqrt{\log T}\big)\qquad(\nu=1,2),
\end{eqnarray}
where $c(A)$ is a positive function of $A$. 

Up to this point, all the analysis is unconditional. To study the
remaining case, in which $\mathscr{L}^A=\eta<q\leq y$, Conrey, Ghosh and Gonek
used the Vaughan identity and the large sieve. Their approach requires
the assumption of GLH (or precisely, an upper bound for averages of
sixth moments of Dirichlet $L$-functions). In the next section, we
shall illustrate how Heath-Brown's generalisation of the Vaughan
identity can be used to obtain unconditionally the following estimate.
\begin{lemma}
We have
\[\mathcal{M}_{\nu,3}\ll_{\varepsilon}
y^{1/3}T^{5/6+\varepsilon}+\eta^{-1/2}T\mathscr{L}^C\qquad(\nu=1,2), \]
for some absolute constant $C>0$ and for any fixed $\ep>0$.
\end{lemma}

We finish the section with the deduction of our theorem. Given
$\vartheta<1/2$, Lemma 1, Lemma 2 and (5)--(7) give 
\begin{eqnarray*}
S_1\sim\frac{T\mathscr{L}^2}{2\pi}
\bigg(\tfrac{1}{2}+\vartheta\int_{0}^{1}P(u)du\bigg) 
\end{eqnarray*}
and
\begin{eqnarray*}
S_2\sim\frac{T\mathscr{L}^3}{2\pi}
\bigg(\tfrac{1}{3}+\vartheta\int_{0}^{1}P(u)du+\vartheta^2
\bigg(\int_{0}^{1}P(u)du\bigg)^2+\tfrac{1}{12\vartheta}\int_{0}^{1}P'(u)^2du\bigg).
\end{eqnarray*}
Choosing $P(x)=-\vartheta x^2+(1+\vartheta)x$ and letting
$\vartheta\rightarrow1/2^{-}$ we obtain 
\begin{equation}\label{8}
S_1\sim\frac{19}{24}\frac{T\mathscr{L}^2}{2\pi}\qquad\textrm{and}\qquad
S_2\sim\frac{57}{64}\frac{T\mathscr{L}^3}{2\pi}. 
\end{equation}
Assuming RH we have $S_2=\sum_{0<\gamma\leq T}|B\zeta'(\rho)|^2$. Note
that this is the only place we need RH. The theorem then follows from
\eqref{1} and \eqref{8}.

\section{Proof of Lemma 2}

We shall prove Lemma 2 for $\mathcal{M}_{2,3}$, the treatment of
$\mathcal{M}_{1,3}$ being similar. 

\subsection{Initial cleaning}

There are problems arising with the condition $d|kq$ in \eqref{9},
since we would like to be able to separate the variables $d$ and $q$.
Indeed it appears that Conrey, Ghosh and Gonek run into difficulties
at this point in deducing [\textbf{\ref{CGG1}}; (7.2)] from
[\textbf{\ref{CGG1}}; (5.15)].  To circumvent such problems we begin
by observing
that the function $b(*)$ given by \eqref{bdef} is supported on
squarefree values.  It follows that we can restrict $k$ and $q$ in
\eqref{9} to be coprime.  Thus the variable $l$ in \eqref{deldef}
will be coprime to $q$.  One then sees that the term $\psi(d/l)$ will
vanish unless $d$ is also coprime to $q$, since $\psi$ is a character 
to modulus $q$. Finally, if $d$ is coprime to $q$, the condition
$d|kq$ reduces to $d|k$.  We therefore conclude that 
\begin{eqnarray}\label{9'}
\mathcal{M}_{2,3}=\sum_{\eta<q\leq y}\sum_{\psi(\textrm{mod}\
  q)}{\!\!\!\!\!\!}^{\textstyle{*}}\
\tau(\overline{\psi})\sum_{k\leq
  y/q}\frac{b(kq)}{kq}\sum_{d|k}\delta(q,kq,d,\psi)\sum_{m\leq
  kqT/2\pi d}a_2(md)\psi(m).
\end{eqnarray}

We divide the summation over $k,q,d$ in \eqref{9'}
into dyadic intervals 
\[K/2<k\leq K,\;\;\; Q/2<q\leq Q,\;\;\; D<d\le 2D,\] 
where
\begin{equation}\label{R2}
Q>\eta=\mathscr{L}^A,\;\;\; D\le K\;\;\;\mbox{and}\;\;\; KQ\leq 4y. 
\end{equation}
Then there will be some such triple $K,Q,D$ for which we have  
\[\mathcal{M}_{2,3}\ll\mathscr{L}^3\sum_{d\sim D}
\sum_{\substack{k\sim K\\d|k}}\sum_{q\sim Q}
|\tau(\overline{\psi})|\frac{|b(kq)|}{kq}|\delta(q,kq,d,\psi)|
\sum_{\psi(\textrm{mod}\ q)}{\!\!\!\!\!\!}^{\textstyle{*}}\
\bigg|\sum_{m\leq kqT/2\pi d}a_2(md)\psi(m)\bigg|.\]
We now note that $|\tau(\psi)|=q^{1/2}$ and
\begin{equation*}
\delta(q,kq,d,\psi)\ll\sum_{l|d}\frac{1}{\varphi(kq/l)}\ll\mathscr{L}
dk^{-1}q^{-1},
\end{equation*}
since $\varphi(n)\gg n/\log\log n$ and $\sigma(n)\ll n\log\log n$. This allows us to write
\begin{eqnarray*}
\mathcal{M}_{2,3}&\ll& K^{-2}Q^{-3/2}D\mathscr{L}^4\sum_{d\sim D}
\sum_{\substack{k\sim K\\d|k}}\sum_{q\sim Q}
\sum_{\psi(\textrm{mod}\ q)}{\!\!\!\!\!\!}^{\textstyle{*}}\
\bigg|\sum_{m\leq kqT/2\pi d}a_2(md)\psi(m)\bigg|\\
&\ll& K^{-2}Q^{-3/2}D\mathscr{L}^4\sum_{d\sim D}
\sum_{\substack{k\sim K\\d|k}}\;S(Q,X,d),
\end{eqnarray*}
where we have defined
\[X=KQT/\pi D\]
and
\[S(Q,X,d)=\sum_{q\sim Q}\;\sum_{\psi(\textrm{mod}\
  q)}{\!\!\!\!\!\!}^{\textstyle{*}}\ \ \max_{M\leq X}\bigg|\sum_{m\leq
  M}a_2(md)\psi(m)\bigg|.\]
Since the number of available values for $k$ is $\ll K/d\ll K/D$
we conclude that
\begin{eqnarray}\label{rep}
\mathcal{M}_{2,3}&\ll&K^{-1}Q^{-3/2}\mathscr{L}^4\sum_{d\sim D}\;
S(Q,X,d).
\end{eqnarray}

The sum $S(Q,X,d)$ would be in a suitable form to apply the
maximal large sieve, if
it involved the square of the innermost sum.  However $X$ is too large
compared with $Q$ for one merely to apply Cauchy's inequality.  Thus
the strategy is to use a generalisation of the Vaughan identity to write
the function $a_2$ as a convolution, thereby enabling us to replace
the innermost sum by a product of two Dirichlet polynomials. Providing
these two polynomials are of suitable lengths a satisfactory estimate
will emerge. The details of our implementation differ from those of
Conrey, Ghosh and Gonek in two important ways.  Firstly, 
by using the identity in
Lemma 3 we produce more flexibility in the choice of lengths for our
Dirichlet polynomials.  Secondly, Conrey, Ghosh and Gonek used
$L(s,\psi)$ where we employ a finite Dirichlet polynomial of the type
$\sum_{h\sim H}h^{-s}\psi(m)$. This is clearly advantageous if $H$ is
small.

There are two inconvenient technical problems which need to be dealt with.
Firstly, since we have $a_2(md)$ rather than merely $a_2(m)$ we
have to handle the dependence on $d$.  Secondly, when we replace
$a_{2}$ by a convolution we need to eliminate the condition $m\le
M$.  We do this in the standard way by using Perron's formula, which
introduces a further variable, and a further averaging, into our analysis.

\subsection{The generalised Vaughan identity}

Heath-Brown's version [\textbf{\ref{H-B1}}] of the Vaughan
identity comes from the following trivial lemma.

\begin{lemma}
For any integer $r\geq1$ we have
\begin{eqnarray}\label{10}
\zeta'(s)/\zeta(s)=\sum_{j=1}^{r}(-1)^{j-1}
\binom{r}{j}\zeta(s)^{j-1}\zeta'(s)M(s)^j+
\big(1-\zeta(s)M(s)\big)^r\zeta'(s)/\zeta(s), 
\end{eqnarray}
where
\begin{equation*}
M(s)=\sum_{n\leq X}\frac{\mu(n)}{n^s}.
\end{equation*}
\end{lemma} 

We apply Lemma 3 to the sum $S(Q,X,d)$, where the coefficients
$a_2(n)$ are defined in \eqref{3}, so that $a_2=-\Lambda*\log*\log*\, b$. We
choose $r=3$, $X=T^{1/2}$, and pick out the relevant coefficients of
$n^{-s}$ with $n=md$. Since 
\begin{equation}\label{R1}
Md\le KQT/\pi\le 4yT/\pi<T^{3/2}
\end{equation}
for large $T$, we see
that the last term on the right hand side of
\eqref{10} makes no contribution. On splitting each range of summation
into dyadic intervals, we find that $a_2(md)$ is a linear
combination of $O(\mathscr{L}^9)$ expressions of the form
$(f_1*\ldots*f_9)(md)$, where the functions $f_i$ are independent of
$m$ and $d$, and are each supported on a dyadic interval
$(N_i/2,N_i]$, say.  For terms in which the function $f_i$ is absent we
set $N_i=1$ and take the corresponding function $f_i$ to be the
identity for the Dirichlet convolution, so that $f_i(1)=1$ and
$f_i(m)=0$ for $m\ge 2$.  Whenever $N_i>1$ we can take 
\[f_1=f_2=f_3=\log,\;\;\; f_4=b,\;\;\; f_5=f_6=1,\;\;\;\mbox{and}
\;\;\; f_7=f_8=f_9=\mu.\]
Moreover
\[N_4\le y\;\;\;\mbox{and}\;\;\; N_7,N_8,N_9\le T^{1/2}.\]
We observe that the numbers $N_i$ run over powers of 2 or, in the case
of $N_4$ over numbers $2^{-h}y$.  Since these are independent of $q$
we can estimate $S(Q,X,d)$ as
\[S(Q,X,d)\ll\sum_{N_i}\sum_{q\sim Q}\;
\sum_{\psi(\textrm{mod}\ q)}{\!\!\!\!\!\!}^{\textstyle{*}}\  \
\max_{M\leq X} \bigg|\sum_{m\leq
  M}(f_1*\ldots*f_9)(md)\psi(m)\bigg|,\]
where the sum over $N_i$ runs through $O(\mathscr{L}^9)$ sets of
values with $\prod N_i\ll Xd$.

To evaluate $(f_1*\ldots*f_9)(md)$ we call on Lemma 3 of
Conrey, Ghosh and Gonek [\textbf{\ref{CGG1}}], which shows that
\[(f_1*\ldots*f_9)(md)=\sum_{d=d_1\ldots d_9}
(g_1*\ldots*g_9)(m),\]
with 
\[g_i(m)=g_i(m;d_1,\ldots,d_i)=\left\{\begin{array}{ll}
f_i(md_i),&\quad \textrm{if }(m,d_1\ldots d_{i-1})=1,\\
0,&\quad \textrm{otherwise.}\\ \end{array}\right.\]
Each $g_i$ is now supported on a dyadic interval $(M_i/2,M_i]$ with
$M_i=N_i/d_i$, so that $\prod M_i\ll X$.

This allows us to estimate $S(Q,X,d)$ as
\[S(Q,X,d)\ll\sum_{N_i}\sum_{d=d_1\ldots d_9}
\sum_{q\sim Q}\;\sum_{\psi(\textrm{mod}\ q)}{\!\!\!\!\!\!}^{\textstyle{*}}\
\ \max_{M\leq X} \bigg|\sum_{m\leq
  M}(g_1*\ldots*g_9)(m)\psi(m)\bigg|.\]

We begin by disposing of the case in which $M_i>yT^{1/2}$ for
some index $i$, which will necessarily be 1, 2, 3, 5 or 6.  
For $B\ll X$ we can use partial summation to show that
\[\sum_{A<h\le B}g_i(h)\psi(h)
\ll\mathscr{L}\max_{Y\le B}\left|\sum_{\substack{h\le Y\\
(h,D_i)=1}}\psi(h)\right|,\]
with $D_i=d_1\ldots d_{i-1}$. Moreover
\[\sum_{\substack{h\le Y\\(h,D_i)=1}}\psi(h)=\sum_{e|D_i}\mu(e)\psi(e)
\sum_{k\le Y/e}\psi(k)\ll_\ep \tau(D_i)q^{1/2}\log q\]
by the P\'olya--Vinogradov inequality. It follows that
\[\sum_{A<h\le B}g_i(h)\psi(h)\ll_\ep Q^{1/2}T^{\ep}\]
for $B\ll X$.

We now write $g$ for
the convolution of the 8 functions $g_j$ with $j\not=i$, so that 
$g$ is supported on integers $n\ll X/M_i$, and $g(n)\ll_\ep T^{\ep}$.  Then
\begin{eqnarray*}
\sum_{m\leq M}(g_1*\ldots*g_9)(m)\psi(m)
&=&\sum_{n}g(n)\psi(n)\sum_{\substack{h\sim M_i\\ h\le M/n}}g_i(h)\psi(h)\\
&\ll_\ep&XM_i^{-1}Q^{1/2}T^{2\ep}\\
&\ll_\ep& XQ^{1/2}y^{-1}T^{-1/2+2\ep}.
\end{eqnarray*}
The contribution to $S(Q,X,d)$ when $M_i>yT^{1/2}$ is therefore
\begin{eqnarray}\label{rev1}
&\ll_\ep& \mathscr{L}^9\tau_9(d)Q^2.XQ^{1/2}y^{-1}T^{-1/2+2\ep}\nonumber\\
&\ll_\ep& KQ^{7/2}D^{-1}y^{-1}T^{1/2+3\ep}
\end{eqnarray}
by \eqref{R1} and \eqref{R2}.

Before handling the remaining terms we must
eliminate the condition $m\le M$, which may be done via
Perron's formula. Let $M_0=M+1/2$ and $\delta=(\log M)^{-1}$, and
take $U>0$.  Then
\[\frac{1}{2\pi i}\int_{\delta-iU}^{\delta+iU}
\left(\frac{M_0}{m}\right)^s\frac{ds}{s}=\left\{\begin{array}{cc} 1
&\mbox{ if } m\le M\\ 0&\mbox{ if } m> M\end{array}\right\}
+O(MU^{-1}).\]
Thus
\begin{eqnarray*}
\lefteqn{S(Q,X,d)}\\
&\ll& 1+\sum_{N_i}\sum_{d=d_1\ldots d_9}\int_{-U}^U\frac{\log X}{1+|t|}
\sum_{q\sim Q}\;
\sum_{\psi(\textrm{mod}\ q)}{\!\!\!\!\!\!}^{\textstyle{*}}\ \ \bigg|\sum_m
(g_1*\ldots*g_9)(m)\psi(m)m^{-\delta-it}\bigg|dt,
\end{eqnarray*}
provided that $U\gg Q^2X^2T^{\delta'}$ for some fixed $\delta'>0$.  The
reader should note here that the sum over $m$ is finite, being supported on
values
\[m=m_1\ldots m_9\le M_1\ldots M_9\ll X.\]
We now choose $U=T^5$, which is more than sufficient when 
$\delta'=1/2$, say.  We proceed to define functions
$h_j(m)=g_j(m)m^{-\delta}$, and set
\[H_j(\psi,t)=\sum_{m\sim M_j}h_j(m)\psi(m)m^{-it},\]
which allows us to conclude that
\begin{equation}\label{E1}
S(Q,X,d)\ll 1+\mathscr{L}^2\sum_{N_i}\sum_{d=d_1\ldots
  d_9}\max_{1\le V\le T^5}V^{-1}T(Q,V),
\end{equation}
where
\[T(Q,V)=T(Q,V;d_1,\ldots,d_9)=
\sum_{q\sim Q}\;
\sum_{\psi(\textrm{mod}\ q)}{\!\!\!\!\!\!}^{\textstyle{*}}\ \ 
\int_{-V}^V|H_1(\psi,t)\ldots H_9(\psi,t)|dt.\]

\subsection{Estimating $T(Q,V)$}

We may now suppose that $M_i\le yT^{1/2}$ for every
index $i$.  Our strategy is to split $\prod H_i(\psi,t)$ into a product
$\mathscr{A}(\psi,t)\mathscr{B}(\psi,t)$ of Dirichlet polynomials of
approximately equal lengths $A$ and $B$ respectively. They will take the form
\[\sum_{m\le A}a_m\psi(m)m^{-it},\;\;\;\mbox{and}\;\;\;
\sum_{m\le B}b_m\psi(m)m^{-it}\]
with coefficients such that 
\begin{equation}\label{R7}
|a_m|,|b_m|\ll
\mathscr{L}^3\tau_9(m). 
\end{equation}
Giving $\mathscr{A}(\psi,t)$ and $\mathscr{B}(\psi,t)$ approximately 
equal lengths will 
optimise our eventual application of the hybrid large sieve. To achieve this 
we introduce a parameter $A_0\ge yT^{1/2}$, to be specified in due course, 
with the aim of making $\max\{A,B\}\ll A_0$. We recall that
\begin{equation}\label{R5}
X= KQT/\pi D.
\end{equation}
Now, if there is a factor $H_i(\psi,t)$ of length $M_i\ge
KQT(DA_0)^{-1}$ we can merely take
$\mathscr{A}(\psi,t)=H_i(\psi,t)$.  We will then have $A=M_i\le
yT^{1/2}\le A_0$.  Moreover, since $A=M_i\ge KQT(DA_0)^{-1}$, the 
corresponding factor $\mathscr{B}(\psi,t)$
will have $B\ll X/A\ll KQT(DA)^{-1}\ll A_0$ as required.  We can therefore
assume that each $H_i(\psi,t)$ has length $M_i\le KQT(DA_0)^{-1}$.

In this remaining case we define $J$ as the 
largest integer for which
$\prod_{j\le J}M_j\le A_0$. We then set
\[\mathscr{A}(\psi,t)=\prod_{j\le J}H_j(\psi,t),\;\;\; \mbox{and}\;\;\;
\mathscr{B}(\psi,t)=\prod_{J<j\le 9}H_j(\psi,t)\]
so that $A\le A_0$.  Moreover our construction implies that
$AM_{J+1}>A_0$, and since we
are assuming that $M_i\le KQT(DA_0)^{-1}$ for every
index $i$ we see that 
\[A\gg\frac{A_0}{KQT(DA_0)^{-1}},\]
whence 
\[B\ll\frac{X}{A}\ll\frac{KQT/D}{A}\ll \left(\frac{KQT}{DA_0}\right)^2.\]
We therefore see that if we set
\begin{equation}\label{R6}
A_0=\max\big\{yT^{1/2}\,,\,(KQT/D)^{2/3}\big\}
\end{equation}
then we can always produce a factorisation with $A,B\ll A_0$.

Having chosen $\mathscr{A}(\psi,t)$ and $\mathscr{B}(\psi,t)$
we proceed to apply Cauchy's inequality to obtain
\[T(Q,V)\le T(\mathscr{A})^{1/2}T(\mathscr{B})^{1/2},\]
where
\[T(\mathscr{A})=\sum_{q\sim Q}\;
\sum_{\psi(\textrm{mod}\ q)}{\!\!\!\!\!\!}^{\textstyle{*}}\ \ 
\int_{-V}^V|\mathscr{A}(\psi,t)|^2dt,\]
and similarly for $T(\mathscr{B})$.  We then use the hybrid large sieve
in the form
\[\sum_{q\sim Q}\;
\sum_{\psi(\textrm{mod}\ q)}{\!\!\!\!\!\!}^{\textstyle{*}}\ \ 
\int_{-V}^V\left|\sum_{m\le H}h_m\right|^2dt\ll (Q^2V+H)\sum|h_m|^2,\]
due to Montgomery [\textbf{\ref{HLM}}; Theorem 7.1].  This produces a bound
\[T(\mathscr{A})
\ll (Q^2V+A)\sum_{m\le A}|a_m|^2\ll (Q^2V+A)A\mathscr{L}^{86},\]
in view of \eqref{R7}, and similarly for $T(\mathscr{B})$.  It
follows that
\begin{eqnarray*}
T(Q,V)&\ll& \big\{(Q^2V+A)A\big\}^{1/2}\big\{(Q^2V+B)B\big\}^{1/2}\mathscr{L}^{86}\\
&\ll&
\big\{Q^2VX^{1/2}+QV^{1/2}X^{1/2}\max\{A^{1/2},B^{1/2}\}+X\big\}\mathscr{L}^{86}\\
&\ll&
\big\{Q^2VX^{1/2}+QV^{1/2}X^{1/2}A_0^{1/2}+X\big\}\mathscr{L}^{86}.
\end{eqnarray*}
Since $V\ge 1$ we now deduce from \eqref{R5} and \eqref{R6} that
\begin{eqnarray}\label{x1}
V^{-1}T(Q,V)&\ll_\ep&
K^{1/2}Q^{5/2}D^{-1/2}T^{1/2+\ep}+K^{1/2}Q^{3/2}D^{-1/2}y^{1/2}T^{3/4+\ep}\nonumber\\
&&\hspace{2cm}+K^{5/6}Q^{11/6}D^{-5/6}T^{5/6+\ep}+
KQD^{-1}T\mathscr{L}^{86}
\end{eqnarray}
for any fixed $\ep>0$, when $\max M_i\le yT^{1/2}$.

\subsection{Deduction of Lemma 2}

Putting the estimate \eqref{x1} into \eqref{E1} and comparing with
\eqref{rev1} we get
\begin{eqnarray*}
S(Q,X,d)&\ll_\ep& KQ^{7/2}D^{-1}y^{-1}T^{1/2+3\ep}+
K^{1/2}Q^{5/2}D^{-1/2}T^{1/2+2\ep}\\
&&\hspace{2cm}+K^{1/2}Q^{3/2}D^{-1/2}y^{1/2}T^{3/4+2\ep}
+K^{5/6}Q^{11/6}D^{-5/6}T^{5/6+2\ep}\\
&&\hspace{3cm}+KQ\tau_9(d)D^{-1}T\mathscr{L}^{97},
\end{eqnarray*}
whence \eqref{R2} and \eqref{rep} yield
\[\mathcal{M}_{2,3}\ll_\ep Q^2y^{-1}T^{1/2+4\ep}+QT^{1/2+3\ep}
+y^{1/2}T^{3/4+3\ep}+Q^{1/3}T^{5/6+3\ep}+
Q^{-1/2}T\mathscr{L}^{109}.\]
Thus since $\eta\ll Q\ll y$ we have
\[
\mathcal{M}_{2,3}\ll_\ep
yT^{1/2+4\ep}+y^{1/2}T^{3/4+3\ep}+y^{1/3}T^{5/6+3\ep}+
\eta^{-1/2}T\mathscr{L}^{109},
\]
and since $y\le T^{1/2}$ the bound required for Lemma 2 follows, on 
re-defining $\ep$.\\

\textbf{Acknowledgement.} We would like to thank Micah Milinovich for a helpful remark.

\end{document}